\documentclass[12pt, twoside, fleqn, a4paper]{article}
\usepackage{lineno}

\usepackage{latexsym}  
\usepackage{amscd}    
\usepackage{theorem}  
\usepackage{pifont}  
\usepackage{mathbbol} 
\usepackage{amsfonts}  
\usepackage{xspace}  
\usepackage{amssymb} 
\usepackage{fancyhdr} 
\usepackage{mathrsfs} 
\usepackage{amsmath} 
\usepackage{graphics} 
\usepackage{graphicx} 
\usepackage{euscript}
\usepackage{psfrag}
\usepackage{empheq} 
\usepackage{mathabx} 
\usepackage[parfill]{parskip}     
\usepackage[colorlinks=true, allcolors=blue]{hyperref}
\usepackage{cancel} 
\usepackage{bigints} 

\title{Extremal Lyapunov exponents in \\ random and simultaneous dynamics}
\author{Thirupathi Perumal\footnote{Indian Institute of Science Education and Research Thiruvananthapuram (IISER-TVM), Maruthamala P.O., \\ Vithura, Kerala, India. PIN 695 551.\ \ email: \texttt{thirupathip23@iisertvm.ac.in}}\ \ and\ \ Shrihari Sridharan\footnote{Indian Institute of Science Education and Research Thiruvananthapuram (IISER-TVM), Maruthamala P.O., \\ Vithura, Kerala, India. PIN 695 551.\ \ email: \texttt{shrihari@iisertvm.ac.in} (Corresponding Author). \\ The second author thanks the support provided by NBHM through a research grant No. \\ 02011/35/2025/NBHM(R.P)/R\&D II/9832}}

\date{\today} 

\DeclareFontFamily{OT1}{pzc}{}
\DeclareFontShape{OT1}{pzc}{m}{it}%
              {<-> s * [0.900] pzcmi7t}{}
\DeclareMathAlphabet{\mathpzc}{OT1}{pzc}%
                                 {m}{it}

{\theorembodyfont{\slshape} \newtheorem{theorem}{Theorem}[section]}
{\theorembodyfont{\slshape} \newtheorem{definition}[theorem]{Definition}}
{\theorembodyfont{\slshape} \newtheorem{lemma}[theorem]{Lemma}}
{\theorembodyfont{\slshape} }
{\theorembodyfont{\slshape} }
{\theorembodyfont{\slshape} } 
{\theorembodyfont{\slshape} } 
{\theorembodyfont{\slshape} }
{\theorembodyfont{\slshape} }

\numberwithin{equation}{section}

\newenvironment{proof}{\paragraph{Proof:}}{\hfill$\bullet$}

\begin{document}

\maketitle

\begin{abstract} 
\noindent 
In this manuscript, we consider finitely many maps, all of which are defined on a compact probability space, with at least one map in the collection having degree strictly bigger than $1$. Working with two different, however related, dynamics generated by this setting, we obtain an expression for the extremal Lyapunov exponents, that characterise the instability of typical orbits, as the limit of the averages of the logarithm of the operator norm of linear cocycles of generic orbits. We obtain this as a consequence to the Kingman's ergodic theorem for a subadditive sequence of measurable functions, which naturally generalises the Birkhoff's ergodic theorem, in each of the cases. 
\end{abstract}

\begin{tabular}{l l} 
{\bf Keywords} & Random dynamical systems, \\ 
& Lyapunov exponents, \\ 
& Subadditive sequence of measurable functions. \\ 
& \\ 
{\bf MSC Subject} & 37H15, 37H12, 37C35, 37A05 \\ 
{\bf Classifications} & \\ 
\end{tabular} 
\bigskip 


\section{Introduction} 
Dynamics is the investigation of evolution of a system in question, over time. Even though studies of such systems are carried out wherein a single map takes charge of the necessary evolution by several authors in many papers, mathematicians have also been interested in systems where more than one map is placed at their disposal for the necessary evolution. It is always interesting to adapt and arrive at well-known and established results available in \emph{single-map-dynamical-system} to the case of \emph{multiple-maps-dynamical-system}. For example, the famous pointwise ergodic theorem due to Birkhoff, as one may find in \cite{pw:2000} says upon satisfying certain natural hypothesis, the local time averages converge to the global space average of the system, almost everywhere in the space, in the case of a single-map-dynamical-system. Recently, the authors in \cite{aswin:2024} proved the Birkhoff's ergodic theorem in the case of multiple-maps-dynamical-systems. Here, the orbits of generic points are allowed to evolve either along a specified branch from the multiple branches available with equal probability (called \emph{random dynamics}) or along the multiple branches with equal probability, resulting in set-valued maps (called \emph{simultaneous dynamics}). Please see Section \ref{prelims} for a precise definition of these terms.  

A natural generalisation of the Birkhoff's ergodic theorem is given by the Kingman's ergodic theorem, wherein one considers a subadditive sequence of measurable observables that takes value in $[- \infty, \infty)$. The theorem due to Kingman, as in \cite{kingman:1968}, then speaks about the convergence of the ratio of this sequence of functions to the time that has elapsed. In a dynamical system where the Kingman ergodic theorem is satisfied, the authors in \cite{furs:1960} define a particular observable function by associating a non-singular matrix to every point in the phase space and obtain an expression for the extremal Lyapunov exponents almost everywhere. What we intend to do in this paper is to exploit the analogue of the Birkhoff's ergodic theorem for random dynamics and simultaneous dynamics, as in \cite{aswin:2024}, generalise the same and obtain the analogues of the Kingman ergodic theorem, as in \cite{kingman:1968} and the Furstenberg-Kesten theorem, as in \cite{furs:1960}, however in our current setting of random and simultaneous dynamics. 

Connoisseurs of dynamical systems may be aware of the recent works by Barrientos and Malicet, where they consider multiple-maps-dynamical-systems, in which the evolution takes place using finitely many $\mathcal{C}^{r}$-diffeomorphisms on a finite volume Riemannian manifold of dimension $d > 1$, in \cite{bm:2020} and Lipschitz maps defined on a compact metric space, in \cite{bm:pp} and obtain an expression for the extremal Lyapunov exponents. In \cite{bm:pp}, they define the annealed Koopman operator on the space of continuous functions defined on $X$, prove it to be quasi-compact and use the same to obtain a generalised Kingman's subadditive ergodic theorem, however for random dynamical systems. Our work is different from \cite{bm:pp}, in the sense that even though the path through which the dynamics evolves is prescribed in the initial part of the paper, as one may infer from Theorems \ref{thm:fkomega} and \ref{kingmanomega}, we also focus on the set-valued simultaneous dynamical systems in the latter part of the paper, as can be seen from Theorems \ref{thm:fkomegafree} and \ref{kingmanomegafree}. Also, the proofs of the results regarding the expression for the extremal Lyapunov exponents in \cite{bm:pp} and the current paper, as well as the proof of the Kingman's ergodic theorem, an essential ingredient in both the papers are totally different. However, it must be mentioned that the work by Barrientos and Malicet, in \cite{bm:pp}, contain far more equivalent conditions for these limits to exist and also has several interesting results pertaining to random dynamics that we have not considered, in this work. Interested readers may find \cite{lq:1995, al:1998} to be a good initial point to learn the basics of multiple-maps-dynamical-systems, while \cite{dk:2016, mv:2014} is a good source on the topic of Lyapunov exponents. 

This manuscript is organised as follows: In Section \ref{prelims}, we explain the basic setting of this paper and explain the necessary background that makes this paper comprehensible. After we recall the definitions of extremal Lyapunov exponents and the theorems due to Furstenberg and Kesten, as written in Theorem \ref{thm:fkonemap} and Kingman, as written in Theorem \ref{kingman}, we state the main theorems of this paper in Section \ref{maintheorems}. The main theorems are written as Theorems \ref{thm:fkomega} and \ref{kingmanomega} for random dynamics and Theorems \ref{thm:fkomegafree} and \ref{kingmanomegafree} for simultaneous dynamics. Once the main theorems are stated, we prove each of them in the subsequent sections, namely Sections \ref{kignmansec}, \ref{proofcontd1} and \ref{proofcontd2}. 

\section{Preliminaries} 
\label{prelims} 

Consider finitely many maps, namely $T_{1}, T_{2}, \cdots, T_{N},\ 1 < N < \infty$, all of which act on $M$, a compact probability space. Let $\mu$ be a non-atomic probability measure supported on $M$, that remains invariant under the action of $T_{j}$ for every $1 \le j \le N$. For example, one may consider the maps $T_{j} (x) = (j + 1) x (\mathrm{mod} 1)$ on $[0, 1)$, all of which preserve the Lebesgue measure therein. The random dynamics that we investigate in this paper, grows as follows on $M$: Any point $x_{0} \in M$ at time $0$ has an equal probability of being at $x_{n} \in M$ at time $n \in \mathbb{Z}_{+}$, where 
\[ x_{n} \in \Big\{ \left( T_{\omega_{n}} \circ \cdots \circ T_{\omega_{1}} \right) (x_{0})\ :\ \omega_{j} \in \left\{ 1, 2, \cdots, N \right\} \Big\}. \] 
It is then a simple observation that we have a possibility of $N^{n}$ points, that one may consider as the $n$-th order forward image of $ x_{0}$, thus rendering a richer dynamical system than dealing with a single map on $M$. As one may note, the choice of the combinatorial $n$-lettered word from the letters $\{ 1, 2, \cdots, N \}$ dictates the position of the evolution of the forward orbit of $x_{0}$. Thus defining, 
\[ \Sigma_{N}^{+}\ =\ \Big\{ \omega = \left( \omega_{1} \omega_{2} \cdots \right)\ :\ \omega_{j} \in \left\{ 1, 2, \cdots, N \right\} \Big\}\ \ \ \ \text{and}\ \ \ \ \pi_{n} \omega\ =\ \left( \omega_{1} \omega_{2} \cdots \omega_{n} \right)\ \ \text{for any}\ n \in \mathbb{Z_{+}}, \] 
one may consider any specific path of any length and study the evolution of the orbit of points in $M$ under the map $T_{\pi_{n} \omega}$ given by $T_{\omega_{n}} \circ \cdots \circ T_{\omega_{1}}$ or consider the various possibilities of the point $x_{n} \in M$ pertaining to every $x_{0} \in M$ and every $T_{\pi_{n} \omega}$, when our focus is on simultaneous dynamics. Finally, we let $n \to \infty$ to learn about the eventual evolution of the orbit of $x_{0}$ in $M$. 

Note, from \cite{poli:1998, pw:2000} that $\Sigma_{N}^{+}$ is a metric space when equipped with the metric $d \left( \omega,\, \tau \right) = \dfrac{1}{2^{\rho \left( \omega,\, \tau \right)}}$ where $\rho \left( \omega,\, \tau \right) = \min \left\{ j \in \mathbb{Z}_{+} : \omega_{j} \ne \tau_{j} \right\}$. The infinite concatenation of the $n$-lettered word $\pi_{n} \omega$ to itself converges to the infinite-lettered word $\omega \in \Sigma_{N}^{+}$, as $n \to \infty$. Further, the open sets here can be written as a countable union of cylinder sets, which are themselves both closed and open. Thus, the cylinder sets form a basis for the topology on $\Sigma_{N}^{+}$. By
\emph{a cylinder set of length $m$}, we mean the set denoted by
\[ \left[ \tau_{1} \cdots \tau_{m} \right]\ \ =\ \ \left\{ \omega \in \Sigma_{N}^{+} : \omega_{1} = \tau_{1};\ \cdots;\ \omega_{m} = \tau_{m} \right\} \] 
where we fix the initial $m$ co-ordinates. We wish to draw the attention of the readers to the position of the cylinder sets, which can occur anywhere in general. However, we necessitate the cylinder sets to be placed at the initial co-ordinates, as defined. Under the topology defined, $\Sigma_{N}^{+}$ is a totally disconnected, compact, perfect metric space on which the shift map $\sigma : \Sigma_{N}^{+} \longrightarrow \Sigma_{N}^{+}$ defined by $\sigma (\omega) = \left( \omega_{2} \omega_{3} \cdots \right)$ is a non-invertible, continuous surjection, that has $N$ local inverse branches, for any $\omega$. These cylinder sets form a semi-algebra that would, in turn generate the Borel sigma-algebra on $\Sigma_{N}^{+}$. The equidistributed Bernoulli measure $\mathfrak{B}$ on a cylinder set of length $m$ is defined as $\displaystyle{\mathfrak{B} \left( \left[ \tau_{1} \cdots \tau_{m} \right] \right) = \frac{1}{N^{m}}}$. 

As explained earlier, our ability to choose any map at random from the collection $\left\{ T_{1}, \cdots, T_{N} \right\}$ to progress in our orbit defines \emph{random dynamics} denoted by $R : \Sigma_{N}^{+} \times M \longrightarrow \Sigma_{N}^{+} \times M$ given by $R \left( \omega,\, x \right) = \left( \sigma \omega,\, T_{\pi_{1} \omega} x \right)$. Considering a measurable map $L : M \longrightarrow {\rm GL}_{d} \left( \mathbb{R} \right)$, we also define the map $\mathfrak{R} : \Sigma_{N}^{+} \times M \times \mathbb{R}^{d} \longrightarrow \Sigma_{N}^{+} \times M \times \mathbb{R}^{d}$ given by $\displaystyle{\mathfrak{R} \left( \omega,\, x,\, v \right) = \left( \sigma \omega,\, T_{\pi_{1} \omega} x,\, L(x) v \right)}$. By virtue of the definitions of $R$ and $\mathfrak{R}$, one can find their $m$-th iterate to satisfy 
\begin{eqnarray} 
\label{iteratesofR}
R^{m} \left( \omega,\, x \right) & = & \left( \sigma^{m} \omega,\, T_{\pi_{m} \omega} x \right) \\ 
\label{iteratesofranR}
\mathfrak{R}^{m} \left( \omega,\, x,\, v \right) & = & \left( \sigma^{m} \omega,\, T_{\pi_{m} \omega} x,\, \prod_{j\, =\, 0}^{m - 1} \left( L \left( T_{\pi_{j} \omega} x \right) \right) v \right), 
\end{eqnarray} 
where we consider $T_{\pi_{0} \omega} (x) \equiv x$, for any $\omega \in \Sigma_{N}^{+}$. Further, since the maps are not assumed to be commutative and matrix multiplication is not commutative, in general, we shall always understand this notation to give a meaningful representation, in the natural order arising from the dynamical system, as defined; meaning 
\[ \prod_{j\, =\, 0}^{m - 1} \left( L \left( T_{\pi_{j} \omega} x \right) \right)\ \ =\ \ \left( L \left( T_{\pi_{m - 1} \omega} x \right) \right) \cdots \left( L \left( T_{\pi_{1} \omega} x \right) \right) \left( L \left( T_{\pi_{0} \omega} x \right) \right). \]

Our objective in this paper is to understand random dynamics that arises due to $\mathfrak{R}$ defined on $\Sigma_{N}^{+} \times M \times \mathbb{R}^{d}$, and use the same to investigate the set-valued dynamics called \emph{simultaneous dynamics} due to $\mathfrak{S}$ that we shall now define on subsets of $M \times \mathbb{R}^{d}$. For any set $E \subseteq M$ and $F \subseteq \mathbb{R}^{d}$, we define 
\begin{eqnarray} 
\label{actionofS} 
\mathfrak{S} \left( E,\, F \right) & = & \bigcup_{\omega_{1}\, =\, 1}^{N}\; \bigcup_{x\, \in\, E}\; \bigcup_{v\, \in\, F}\; \left( T_{\omega_{1}} x,\, \left( L (x) \right) v \right) \nonumber \\ 
& = & \bigcup_{ \left\{ \pi_{1} \omega\; :\; \omega\, \in\, \Sigma_{N}^{+} \right\}}\; \bigcup_{x\, \in\, E}\; \bigcup_{v\, \in\, F}\; \left( T_{\pi_{1} \omega} x,\, \left( L (x) \right) v \right), 
\end{eqnarray} 
so that the action of $\mathfrak{S}^{m}$ on the ordered pair $\left( E,\, F \right)$ can be recursively defined for any $m \in \mathbb{Z}_{+}$ and given by 
\begin{eqnarray*} 
\mathfrak{S}^{m} \left( E,\, F \right) & = & \mathfrak{S} \left( \mathfrak{S}^{m - 1} \left( E,\, F \right) \right) \\ 
& = & \bigcup_{\left\{ \pi_{m} \omega\; :\; \omega\, \in\, \Sigma_{N}^{+} \right\}}\; \bigcup_{x\, \in\, E}\; \bigcup_{v\, \in\, F}\; \left( T_{\pi_{m} \omega} x,\, \prod_{j\, =\, 0}^{m - 1} \left( L \circ T_{\pi_{j} \omega} (x) \right) v \right). 
\end{eqnarray*} 
An easy combinatorial calculation reveals that the set of all $m$-lettered words that one could form with the $N$ letters $\{ 1, 2, \cdots, N \}$ is only $N^{m}$, and hence, the leftmost union in the $m$-th iterate of $\mathfrak{S}$ only poses the problem of taking a finite union over these $N^{m}$ many words and not actually on the uncountable set of all $\omega \in \Sigma_{N}^{+}$. To highlight the same, we shall also define the set-valued map $S : M \longrightarrow 2^{M}_{{\rm finite}}$ where $2^{M}_{{\rm finite}}$ refers to the collection of all finite subsets of $M$. To suit our purposes, we shall define 
\begin{equation} 
\label{actionofsimpleS} 
S^{m} (x)\ \ =\ \ \bigcup_{\left\{ \pi_{m} \omega\; :\; \omega\, \in\, \Sigma_{N}^{+} \right\}}\; T_{\pi_{m} \omega} x\ \ \ \text{for any}\ m \in \mathbb{Z}_{+}. 
\end{equation} 

\section{Main theorems} 
\label{maintheorems} 

In \cite{furs:1960, mv:2014}, the authors consider a single map, say $T$ to act on $M$ and define what are called \emph{extremal Lyapunov exponents} with respect to the evolution of dynamics using this map $T$. To state the same, we need a few notations. For any $A \in {\rm GL}_{d} \left( \mathbb{R} \right)$, consider the operator norm on ${\rm GL}_{d} \left( \mathbb{R} \right)$ given by 
\[ \| A \|_{{\rm op}}\ \ =\ \ \sup_{v\, \in\, \mathbb{R}^{d}\; :\; \|v\|\, =\, 1} \|Av\|,\ \ \ \text{where}\ \ \| \cdot \|\ \text{denotes the Euclidean norm on}\ \mathbb{R}^{d}. \] 
For the map $L : M \longrightarrow {\rm GL}_{d} \left( \mathbb{R} \right)$, let $\left( L(x) \right)^{\pm 1}$ denote the matrix $L(x)$ or its inverse, as appropriate. Assume that $\log^{+} \| \left( L(x) \right)^{\pm 1} \|_{{\rm op}}$ is an integrable function with respect to the measure $\mu$, where $\log^{+} \| A \|_{{\rm op}} = \max \left\{ 0,\; \log \| A \|_{{\rm op}} \right\}$. 

As the authors indicate in \cite{aev:pp2023}, one may view the linear cocycle determined by $L$, as a generalisation of random products of matrices from ${\rm GL}_{d} \left( \mathbb{R} \right)$. The authors, in \cite{aev:pp2023}, have investigated the continuity of Lyapunov exponents, associated with such random products. The specific case of $d = 2$ was dealt by the authors, in \cite{tv:2020} while the authors, in \cite{bv:2017} considered analogous problems for random products of matrices from ${\rm GL}_{d} \left( \mathbb{C} \right)$. 

\begin{theorem}[Furstenberg-Kesten Theorem]\cite{furs:1960, mv:2014} 
\label{thm:fkonemap}
Let $T$ be a measure preserving map defined on a compact probability measure space $(M, \mu)$. Suppose $L : M \longrightarrow {\rm GL}_{d} \left( \mathbb{R} \right)$ is a measurable map such that $\log^{+} \left\| L(x)^{\pm 1} \right\|_{{\rm op}} \in \mathscr{L}^{1} (\mu)$. Then, the extremal Lyapunov exponents defined by 
\begin{eqnarray*} 
\lambda_{+} (x) & = & \lim_{m\, \to\, \infty} \frac{1}{m} \log \left\| L \left( T^{m - 1} x \right) \cdots L \left( T x \right) L(x) \right\|_{{\rm op}}\ \ \ \ \text{and} \\  
\lambda_{-} (x) & = & \lim_{m\, \to\, \infty} \frac{1}{m} \log \left\| \left( L \left( T^{m - 1} x \right) \cdots L \left( T x \right) L(x) \right)^{-1} \right\|_{{\rm op}}^{-1}, 
\end{eqnarray*} 
exist for $\mu$-almost every $x \in M$. The functions $\lambda_{+}$ and $\lambda_{-}$ are $T$-invariant, meaning $\lambda_{+} \circ T = \lambda_{+}$ and $\lambda_{-} \circ T = \lambda_{-}$, for $\mu$-almost every $x \in M$. Also the functions $\lambda_{+}$ and $\lambda_{-}$ are integrable with respect to $\mu$. Moreover, 
\begin{eqnarray*} 
\int \lambda_{+} \mathrm{d}\mu & = & \lim_{m\, \to\, \infty} \frac{1}{m} \int \log \left\| L \left( T^{m - 1} x \right) \cdots L \left( T x \right) L(x) \right\|_{{\rm op}} \mathrm{d}\mu\ \ \text{and} \\ 
\int \lambda_{-} \mathrm{d}\mu & = & \lim_{m\, \to\, \infty} \frac{1}{m} \int \log \left\| \left( L \left( T^{m - 1} x \right) \cdots L \left( T x \right) L(x) \right)^{-1} \right\|_{{\rm op}}^{-1} \mathrm{d}\mu. 
\end{eqnarray*} 
\end{theorem} 

In order to prove Theorem \ref{thm:fkonemap}, the author in \cite{mv:2014} considers the notion of subadditivity of a sequence of measurable real-valued functions, say $f_{n}$, which we describe below. 

\begin{definition} 
A sequence of measurable functions, say $f_{n} : M \longrightarrow [ - \infty, \infty)$ for $n \ge 1$ is said to be \emph{subadditive with respect to the map $T$} if $f_{n + p} \le f_{n} + \left( f_{p} \circ T^{n} \right)$ for all $n, p \ge 1$. 
\end{definition} 

It is worth mentioning that the authors, in \cite{as:1978} consider a Markov operator $T$ defined on $\mathscr{L}^{1} (M)$, a sequence of superadditive positive integrable functions, say $f_{n}$ with respect to $T$ and derive a version of the ergodic theorem. Here, one should understand superadditivity with the inequality reversed in the definition of subadditivity. For a subadditive sequence of functions, we now state the Kingman's ergodic theorem, as may be found in \cite{kingman:1968, mv:2014}. 

\begin{theorem}[Kingman's ergodic theorem]\cite{kingman:1968, mv:2014} 
\label{kingman}
Let $f_{n} : M \longrightarrow [ - \infty, \infty)$ for $n \ge 1$ be a subadditive sequence of measurable functions, with respect to $T$ such that $f_{1}^{+} = \max \left\{ f_{1}, 0 \right\} \in \mathscr{L}^{1} (\mu)$. Then, the sequence $\displaystyle{\left\{\dfrac{f_{n}}{n}\right\}_{n\, \ge\, 1}}$ converges $\mu$-almost everywhere to some measurable function, say $f : M \longrightarrow [-\infty, \infty)$, satisfying $f \circ T = f,\ \mu$-almost everywhere. Moreover, the positive part of the limit function is integrable and 
\[ \int f \mathrm{d}\mu\ \ =\ \ \lim_{n\, \to\, \infty} \frac{1}{n} \int f_{n} \mathrm{d}\mu\ \ =\ \ \inf_{n\, \ge\, 1} \frac{1}{n} \int f_{n} \mathrm{d}\mu\ \ \in\ \ [-\infty, \infty). \] 
\end{theorem} 

We focus on generalising the ambit of the above two theorems, in our setting of random dynamics and simultaneous dynamics. In other words, we focus on obtaining analogues of Theorems \ref{thm:fkonemap} and \ref{kingman} for the dynamical system that evolves, along the path determined by some $\omega \in \Sigma_{N}^{+}$ or simultaneously on appropriate subsets, as defined in Section \ref{prelims}. 

Making use of the notation mentioned in Equation \eqref{iteratesofranR}, we now define for any $\omega \in \Sigma_{N}^{+}$ and $x \in M$, the quantities 
\begin{eqnarray} 
\label{lambdaplmiomega} 
\lambda_{+} (\omega,\, x) & = & \lim_{m\, \to\, \infty} \frac{1}{m} \log \left\| \prod_{j\, =\, 0}^{m - 1} \left( L \left( T_{\pi_{j} \omega} x \right) \right) \right\|_{{\rm op}}\ \ \ \ \text{and} \nonumber \\ 
\lambda_{-} (\omega,\, x) & = & \lim_{m\, \to\, \infty} \frac{1}{m} \log \left\|\left( \prod_{j\, =\, 0}^{m - 1}  L \left( T_{\pi_{j} \omega} x \right) \right)^{-1} \right\|_{{\rm op}}^{-1}, 
\end{eqnarray} 
provided the limits exist. Observe that the definitions of $\lambda_{+}$ and $\lambda_{-}$, as written in Equation \eqref{lambdaplmiomega} makes sense, by virtue of the evolution of the dynamics of $\mathfrak{R}$. For the dynamics defined by $R$ as given in Equation \eqref{iteratesofR}, a measurable function $\phi : \Sigma_{N}^{+} \times M \longrightarrow [ - \infty, \infty)$ is said to be $R$-invariant, if $\phi \circ R = \phi$ for $\left( \mathfrak{B} \times \mu \right)$-almost every $(\omega,\, x) \in \Sigma_{N}^{+} \times M$. With these definitions, we now state the first main result of this paper, the analogue of the Furstenberg-Kesten theorem for random dynamics. 

\begin{theorem} 
\label{thm:fkomega}
Consider the random dynamics $\mathfrak{R}$, as defined in Equation \eqref{iteratesofranR}, generated by the finitely many measure preserving maps $\left\{ T_{1}, T_{2}, \cdots, T_{N} \right\},\ 1 \le N < \infty$ acting on a compact probability measure space $(M, \mu)$. Suppose $L : M \longrightarrow {\rm GL}_{d} \left( \mathbb{R} \right)$ is a measurable map such that $\log^{+} \left\| L(x)^{\pm 1} \right\|_{{\rm op}} \in \mathscr{L}^{1} (\mu)$. Then, the limits $\lambda_{+}$ and $\lambda_{-}$, as defined in Equation \eqref{lambdaplmiomega} exist for $\left( \mathfrak{B} \times \mu \right)$-almost every $\left( \omega,\, x \right) \in \Sigma_{N}^{+} \times M$. The functions $\lambda_{+}$ and $\lambda_{-}$ are $R$-invariant and are integrable with respect to $\left( \mathfrak{B} \times \mu \right)$, where the map $R$ is as defined in Equation \eqref{iteratesofR}. Moreover, 
\begin{eqnarray*} 
\int \lambda_{+} \mathrm{d} \left( \mathfrak{B} \times \mu \right) & = & \lim_{m\, \to\, \infty} \frac{1}{m} \bigintsss \log \left\| \prod_{j\, =\, 0}^{m - 1} \left( L \left( T_{\pi_{j} \omega} x \right) \right) \right\|_{{\rm op}} \mathrm{d} \left( \mathfrak{B} \times \mu \right)\ \ \text{and} \\ 
\int \lambda_{-} \mathrm{d} \left( \mathfrak{B} \times \mu \right) & = & \lim_{m\, \to\, \infty} \frac{1}{m} \bigintss \log \left\| \left(\prod_{j\, =\, 0}^{m - 1}  L \left( T_{\pi_{j} \omega} x \right) \right)^{-1} \right\|_{{\rm op}}^{-1} \mathrm{d} \left( \mathfrak{B} \times \mu \right). 
\end{eqnarray*} 
\end{theorem} 

We call the functions $\lambda_{+}$ and $\lambda_{-}$ as the \emph{extremal Lyapunov exponents with respect to the random dynamics} defined by $\mathfrak{R}$. In the course of proving this theorem, we will come across a subadditive sequence of measurable functions. We now define the same. 

\begin{definition} 
\label{subadomega}
A sequence of measurable functions $\phi_{n} : \Sigma_{N}^{+} \times M \longrightarrow [ - \infty, \infty)$ is said to be \emph{subadditive with respect to $R$} if for every $\left( \omega,\, x \right) \in \Sigma_{N}^{+} \times M$, we have 
\[ \phi_{n + p} \left( \omega,\, x \right)\ \ \le\ \ \phi_{n} \left(\omega,\, x \right) + \phi_{p} \left( R^{n} \left( \omega,\, x \right) \right)\ \ =\ \ \phi_{n} \left(\omega,\, x \right) + \phi_{p} \left( \sigma^{n} \omega,\, T_{\pi_{n} \omega} x \right)\ \ \ \forall n, p \ge 1. \] 
In case the above inequality holds with the equal sign for every $\left( \omega,\, x \right) \in \Sigma_{N}^{+} \times M$, we say the sequence is \emph{additive with respect to $R$}. 
\end{definition} 

We now state our next main theorem, the analogue of the Kingman's ergodic theorem for random dynamics. 

\begin{theorem} 
\label{kingmanomega} 
Consider the random dynamics $\mathfrak{R}$, as defined in Equation \eqref{iteratesofranR}, generated by the finitely many measure preserving maps $\left\{ T_{1}, T_{2}, \cdots, T_{N} \right\},\ 1 \le N < \infty$ acting on a compact probability measure space $(M, \mu)$. Let $\left\{ \phi_{n} \right\}_{n\, \ge\, 1};\ \phi_{n} : \Sigma_{N}^{+} \times M \longrightarrow [ - \infty, \infty)$ be a sequence of measurable functions that is subadditive with respect to the map $R$, as written in Definition \ref{subadomega}. Assume that the positive part of $\phi_{1}$, denoted by $\phi_{1}^{+} \in \mathscr{L}^{1} \left( \mathfrak{B} \times \mu \right)$. Then, the sequence $\displaystyle{\left\{\dfrac{\phi_{n}}{n}\right\}_{n\, \ge\, 1}}$ converges $\left( \mathfrak{B} \times \mu \right)$-almost everywhere to some measurable, $R$-invariant function $\phi : \Sigma_{N}^{+} \times M \longrightarrow [-\infty, \infty)$. Moreover, the positive part of the limit function is integrable and 
\[ \int \phi \mathrm{d} \left( \mathfrak{B} \times \mu \right)\ \ =\ \ \lim_{n\, \to\, \infty} \frac{1}{n} \int \phi_{n} \mathrm{d} \left( \mathfrak{B} \times \mu \right)\ \ =\ \ \inf_{n\, \ge\, 1} \frac{1}{n} \int \phi_{n} \mathrm{d} \left( \mathfrak{B} \times \mu \right)\ \ \in\ \ [-\infty, \infty). \] 
\end{theorem} 

Note that one may obtain Theorem \ref{thm:fkonemap} and Theorem \ref{kingman} as corollaries to Theorem \ref{thm:fkomega} and Theorem \ref{kingmanomega} respectively by either considering the case $N = 1$ or considering the infinite lettered word $\omega$ for which $\omega_{i} = \omega_{j}$ for all $i, j \ge 1$. One other way to consider Theorem \ref{thm:fkonemap} and Theorem \ref{kingman} as corollaries to Theorem \ref{thm:fkomega} and Theorem \ref{kingmanomega} respectively is to observe the similaroty between the map $R$ acting on the compact product space $\Sigma_{N}^{+} \times M$ and the lone map $T$ acting on $M$. 

Finally, we aim to write the analogues of Theorem \ref{thm:fkomega} and Theorem \ref{kingmanomega} in the case of simultaneous (set-valued) dynamics, that should thoroughly exploit the richness in the considered dynamical system. For any finite set $E \subsetneq M$, we define the quantities 
\begin{eqnarray} 
\label{lambdaplusomegafree} 
\Lambda_{+} (E) & = & \lim_{n\, \to\, \infty}\; \frac{1}{n}\, \frac{1}{N^{n}}\, \frac{1}{\# E}\, \sum_{\left\{ \pi_{n} \omega\; :\; \omega\, \in\, \Sigma_{N}^{+} \right\}}\, \sum_{x\, \in\, E}\, \log \left\| \prod_{j\, =\, 0}^{n - 1} L \left( T_{\pi_{j} \omega} (x) \right) \right\|_{{\rm op}}\ \ \ \ \text{and} \nonumber \\ 
\Lambda_{-} (E) & = & \lim_{n\, \to\, \infty}\; \frac{1}{n}\, \frac{1}{N^{n}}\, \frac{1}{\# E}\, \sum_{\left\{ \pi_{n} \omega\; :\; \omega\, \in\, \Sigma_{N}^{+} \right\}}\, \sum_{x\, \in\, E}\, \log \left\| \left( \prod_{j\, =\, 0}^{n - 1} L \left( T_{\pi_{j} \omega} (x) \right) \right)^{-1} \right\|_{{\rm op}}^{-1}, 
\end{eqnarray} 
provided the limits exist. Observe that the definitions of $\Lambda_{+}$ and $\Lambda_{-}$ have now included the totality of all branches of all points in the finite set $E$ through which the dynamics evolves. We now state the analogue of the Furstenberg-Kesten theorem in this setting. 

\begin{theorem} 
\label{thm:fkomegafree} 
Consider the simultaneous dynamics $\mathfrak{S}$, as defined in Equation \eqref{actionofS}, generated by the finitely many measure preserving maps $\left\{ T_{1}, T_{2}, \cdots, T_{N} \right\},\ 1 \le N < \infty$ acting on a compact probability measure space $(M, \mu)$. Suppose $L : M \longrightarrow {\rm GL}_{d} \left( \mathbb{R} \right)$ is a measurable map such that $\log^{+} \left\| L(x)^{\pm 1} \right\|_{{\rm op}} \in \mathscr{L}^{1} (\mu)$. Then, for $\mu$-almost every $x \in M$, defining $E_{0} \equiv \{ x \}$ and $\displaystyle{E_{m} = \bigcup\limits_{\left\{ \pi_{m} \omega\; :\; \omega\, \in\, \Sigma_{N}^{+} \right\}} T_{\pi_{m} \omega} \left( x \right)}$ for every $m \in \mathbb{Z}_{+}$, the limits $\Lambda_{+}$ and $\Lambda_{-}$ exist. Further, the functions $\Lambda_{\pm}$ satisfy $\Lambda_{\pm} \left( E_{0} \right) = \Lambda_{\pm} \left( E_{l} \right)$ for any $l \in \mathbb{Z}_{+}$. Since $E_{0} = \{x\}$ is a singleton, considering $\Lambda_{\pm} \left( \left\{ x \right\} \right) = \Lambda_{\pm} (x)$, we also obtain the functions $\Lambda_{\pm}$ to be integrable with respect to $\mu$. Moreover, 
\begin{eqnarray*} 
\int \Lambda_{+} \mathrm{d}\mu & = & \lim_{n\, \to\, \infty} \frac{1}{n} \frac{1}{N^{n}} \bigintsss \sum_{\left\{ \pi_{n} \omega\; :\; \omega\, \in\, \Sigma_{N}^{+} \right\}} \log \left\| \prod_{j\, =\, 0}^{n - 1} L \left( T_{\pi_{j} \omega} (x) \right) \right\|_{{\rm op}} \mathrm{d}\mu\ \ \text{and} \\ 
\int \Lambda_{-} \mathrm{d}\mu & = & \lim_{n\, \to\, \infty} \frac{1}{n} \frac{1}{N^{n}} \bigintss \sum_{\left\{ \pi_{n} \omega\; :\; \omega\, \in\, \Sigma_{N}^{+} \right\}} \log \left\| \left( \prod_{j\, =\, 0}^{n - 1} L \left( T_{\pi_{j} \omega} (x) \right) \right)^{-1} \right\|_{{\rm op}}^{-1} \mathrm{d}\mu. 
\end{eqnarray*} 
\end{theorem} 

The functions $\Lambda_{+}$ and $\Lambda_{-}$ are called the \emph{extremal Lyapunov exponents for simultaneous dynamics} due to the finite collection of maps, $\left\{ T_{1}, T_{2}, \cdots, T_{N} \right\}$, where as mentioned in Theorem \ref{thm:fkomegafree}, one must understand $\Lambda_{\pm} (x) \equiv \Lambda_{\pm} \left( \left\{ x \right\} \right)$. In the sequel, we shall refer to $\Lambda_{\pm}$ as functions on $M$, with this convention, even if we do not mention the same. 

Our next job at hand is to generalise the concept of invariance and subadditivity for a sequence of real-valued functions where the system now evolves with respect to simultaneous dynamics, as given by the map $\mathfrak{S}$. In order to define subadditivity, consider a function $\Phi : M \longrightarrow [ - \infty, \infty)$ and define for every $m \in \mathbb{Z}_{+}$, 
\begin{equation} 
\label{composition} 
\left( \Phi \circ S^{m} \right) (x)\ \ =\ \ \sum_{\left\{ \pi_{m} \omega\; :\; \omega\, \in\, \Sigma_{N}^{+} \right\}} \Phi \circ T_{\pi_{m} \omega} (x), 
\end{equation} 
where the map $S$ is as defined in Equation \eqref{actionofsimpleS}. Then, we have the following definition. 

\begin{definition} 
\label{ksubadd} 
A sequence of measurable functions $\left\{ \Phi_{n} \right\}_{n\, \ge\, 1};\ \Phi_{n} : M \longrightarrow [ - \infty, \infty)$ is said to be \emph{subadditive with respect to $S$} if for every $x \in M$, we have 
\[ \Phi_{n + p} (x)\ \ \le\ \ \Phi_{n} (x) + \left( \Phi_{p} \circ S^{n} \right) (x),\ \ \forall n, p \ge 1. \] 
\end{definition} 

As one may observe, every point in $M$ has $N$ images and thus, it is only natural to include a factor $N$ when we define invariance under the action of $S$. 

\begin{definition} 
A measurable map $\Phi : M \longrightarrow [ - \infty, \infty)$ is said to be \emph{$S$-invariant} if $\left( \Phi \circ S \right) = N \Phi,\ \mu$-almost everywhere on $M$. 
\end{definition} 

With these definitions, we now state the analogue of the Kingman's ergodic theorem, now in the setting of simultaneous dynamics. 

\begin{theorem} 
\label{kingmanomegafree} 
Consider the simultaneous dynamics $\mathfrak{S}$, as defined in Equation \eqref{actionofS}, generated by the finitely many measure preserving maps $\left\{ T_{1}, T_{2}, \cdots, T_{N} \right\},\ 1 \le N < \infty$ acting on a compact probability measure space $(M, \mu)$. Let $\left\{ \Phi_{n} \right\}_{n\, \ge\, 1};\ \Phi_{n} : M \longrightarrow [ - \infty, \infty)$ be a sequence of measurable functions that is subadditive with respect to $S$, as written in Definition \ref{ksubadd}. Assume that the positive part of $\Phi_{1}$, denoted by $\Phi_{1}^{+} \in \mathscr{L}^{1} \left( \mu \right)$. Then, the sequence $\displaystyle{\left\{\dfrac{\Phi_{n}}{n N^{n}}\right\}_{n\, \ge\, 1}}$ converges to some measurable, $S$-invariant function $\Phi : M \longrightarrow [-\infty, \infty),\ \mu$-almost everywhere. Moreover, the positive part of the limit function is integrable and 
\[ \int \Phi \mathrm{d}\mu\ \ =\ \ \lim_{n\, \to\, \infty} \frac{1}{n} \frac{1}{N^{n}} \int \Phi_{n} \mathrm{d}\mu\ \ =\ \ \inf_{n\, \ge\, 1} \frac{1}{n} \frac{1}{N^{n}} \int \Phi_{n} \mathrm{d}\mu\ \ \in\ \ [-\infty, \infty). \] 
\end{theorem} 

\section{Proof of Theorem \ref{kingmanomega}} 
\label{kignmansec} 

In this section, we focus on proving Theorem \ref{kingmanomega}, which indeed is the Kingman's ergodic theorem, if one were to consider the map $R$ defined on $\Sigma_{N}^{+} \times M$ in place of $T$ defined on $M$. However, we will present a complete proof of Theorem \ref{kingmanomega}, since it will be helpful in the adaptation of the same when we discuss simultaneous dynamics, later. 

\begin{proof}[of Theorem \ref{kingmanomega}] 
Given that $\left\{ \phi_{n} \right\}_{n\, \ge\, 1}$ is a subadditive sequence of measurable functions with respect to $R$, it is an easy observation from Definition \ref{subadomega} that 
\[ \phi_{n} (\omega,\, x)\ \ \le\ \ \phi_{1}(\omega,\, x) + \phi_{1}(\sigma \omega,\, T_{\pi_{1} \omega} x) + \cdots + \phi_{1} (\sigma^{n - 1}\omega,\, T_{\pi_{n - 1} \omega} x). \] 
Hence, $\displaystyle{\int \phi_{n} \mathrm{d} \left( \mathfrak{B} \times \mu \right) \le n \int \phi_{1} \mathrm{d} \left( \mathfrak{B} \times \mu \right)}$. Since $\phi_{1}^{+} \in \mathscr{L}^{1} \left( \mathfrak{B} \times \mu \right)$, we have that the integral $\displaystyle{\int \phi_{n} \mathrm{d} \left( \mathfrak{B} \times \mu \right)}$ exists, for every $n \ge 1$. We now state a lemma from elementary analysis and make use of the same to obtain the concluding assertions of the theorem. 

\begin{lemma}[Fekete's Lemma \cite{coor:2015, mv:2014}] 
\label{lim=inf} 
Suppose $\left\{ a_{k} \right\}$ is a subadditive sequence of numbers in $[- \infty, \infty)$, {\it i.e.}, $a_{n + p} \le a_{n} + a_{p}$ for all $n, p \ge 1$. Then 
\[ \lim_{k\, \to\, \infty} \frac{a_{k}}{k}\ \ =\ \ \inf_{k\, \ge\, 1} \frac{a_{k}}{k}. \] 
\end{lemma} 

Thus, defining $\displaystyle{A_{n} = \int \phi_{n} \mathrm{d}\left( \mathfrak{B} \times \mu \right)}$, we observe that $\left\{ A_{n} \right\}_{n\, \ge\, 1}$ is a subadditive sequence of real numbers and hence from Lemma \ref{lim=inf}, we obtain 
\[ A\ \ =\ \ \lim_{n\, \to\, \infty} \frac{1}{n} A_{n}\ =\ \ \lim_{n\, \to\, \infty} \frac{1}{n} \int \phi_{n} \mathrm{d}\left( \mathfrak{B} \times \mu \right)\ \ =\ \ \inf_{n\, \ge\, 1} \frac{1}{n} \int \phi_{n} \mathrm{d}\left( \mathfrak{B} \times \mu \right). \] 

We now consider 
\[ F \left( \omega,\, x \right)\ =\ \liminf_{n\, \to\, \infty} \frac{\phi_{n}}{n} \left( \omega,\, x \right)\ \ \ \text{and}\ \ \ G \left( \omega,\, x \right)\ =\ \limsup_{n\, \to\, \infty} \frac{\phi_{n}}{n} \left( \omega,\, x \right). \] 

Observe that $F \left( \omega,\, x \right) = F \left( \sigma^{m} \omega,\, T_{\pi_{m} \omega} x \right)$ and $G \left( \omega,\, x \right) = G \left( \sigma^{m} \omega,\, T_{\pi_{m} \omega} x \right)$ for any $m \in \mathbb{Z}_{+}$. Also, note that by our definition of $A, F$ and $G$, we have $\displaystyle{\int F \mathrm{d}\left( \mathfrak{B} \times \mu \right) \le A \le \int G \mathrm{d}\left( \mathfrak{B} \times \mu \right)}$. Suppose we also have $\displaystyle{\int G \mathrm{d}\left( \mathfrak{B} \times \mu \right) \le A\le \int F \mathrm{d}\left( \mathfrak{B} \times \mu \right)}$, then the proof of Theorem \ref{kingmanomega} is complete. Towards that end, we initially focus on the sequence of functions $\phi_{n}$ for which $\phi_{n} \left( \omega,\, x \right) > - \infty$ and $F \left( \omega,\, x \right) > - \infty$ for every $\left( \omega,\, x \right) \in \Sigma_{N}^{+} \times M$ and obtain an upper bound for $\displaystyle{\int F \mathrm{d}\left( \mathfrak{B} \times \mu \right)}$. 

Suppose $\dfrac{\phi_{n} \left( \omega,\, x \right)}{n}$ is uniformly bounded below, {\it i.e.}, $\dfrac{\phi_{n} \left( \omega,\, x \right)}{n} \ge - \alpha > - \infty$ for every $ \left( \omega,\, x \right) \in \Sigma_{N}^{+} \times M$ for some $\alpha \in \mathbb{R}_{+}$. Then, $F \left( \omega,\, x \right) \ge - \alpha$. Applying Fatou's lemma, as one may find in \cite{halmos:1970}, to the sequence of non-negative functions $\left\{ \dfrac{\phi_{n} \left( \omega,\, x \right)}{n} + \alpha \right\}_{n\, \ge\, 1}$, we have 
\begin{eqnarray}
\label{oneway} 
\int \left( F + \alpha \right) \mathrm{d}\left( \mathfrak{B} \times \mu \right) & \le & \liminf_{n\, \to\, \infty} \int \left( \dfrac{\phi_{n} \left( \omega,\, x \right)}{n} + \alpha \right) \mathrm{d}\left( \mathfrak{B} \times \mu \right) \nonumber \\ 
& \le & \lim_{n\, \to\, \infty} \int \left( \dfrac{\phi_{n} \left( \omega,\, x \right)}{n} + \alpha \right) \mathrm{d}\left( \mathfrak{B} \times \mu \right) \nonumber \\ 
& = & A + \alpha. 
\end{eqnarray} 

In order to find a lower bound for $\displaystyle{\int F \mathrm{d}\left( \mathfrak{B} \times \mu \right)}$, we fix $\epsilon > 0$ and define a nested sequence of sets, say $\left\{ \mathcal{H}_{k} \right\}_{k\, \in\, \mathbb{Z}_{+}}$ exhausting $\Sigma_{N}^{+} \times M$, given by 
\begin{equation} 
\label{hkw}
\mathcal{H}_{k}\ \ =\ \ \left\{ \left( \omega,\, x \right) \in \Sigma_{N}^{+} \times M\ :\ \phi_{j}\left( \omega,\, x \right) \le j \left( F \left( \omega,\, x \right) + \epsilon \right)\ \ \text{for some}\ 1 \le j \le k \right\}. 
\end{equation} 
Then, defining $\mathcal{F}_{k} : \Sigma_{N}^{+} \times M \longrightarrow \mathbb{R}$ by 
\begin{equation} 
\label{fkw}
\mathcal{F}_{k} \left( \omega,\, x \right)\ \ =\ \ \begin{cases} F \left( \omega,\, x \right) + \epsilon & \text{if}\ \left( \omega,\, x \right) \in \mathcal{H}_{k} \\ 
\phi_{1} \left( \omega,\, x \right) & \text{if}\ \left( \omega,\, x \right)  \in \Sigma_{N}^{+} \times M \setminus \mathcal{H}_{k}, \end{cases} 
\end{equation}
we observe that $\mathcal{F}_{k} \left( \omega,\, x \right) \to F \left( \omega,\, x \right) + \epsilon$ for every $\left( \omega,\, x \right) \in \Sigma_{N}^{+} \times M$. Owing to the definition of $\mathcal{H}_{k}$, we also note that for any $\left( \omega,\, x \right) \in \mathcal{H}_{k + 1} \setminus \mathcal{H}_{k}$, we have 
\[ \mathcal{F}_{k} \left( \omega,\, x \right)\ =\ \phi_{1} \left( \omega,\, x \right)\ >\ F \left( \omega,\, x \right) + \epsilon\ \ \ \text{whereas,}\ \ \ \mathcal{F}_{k + 1} \left( \omega,\, x \right)\ = F \left( \omega,\, x \right) + \epsilon. \] 
Thus, using the monotone convergence theorem, as in \cite{halmos:1970, rudin:1987}, it is easy to infer that 
\[ \int \mathcal{F}_{k} \mathrm{d}\left( \mathfrak{B} \times \mu \right)\ \ \to\ \ \int F \left( \omega,\, x \right) \mathrm{d}\left( \mathfrak{B} \times \mu \right) + \epsilon. \] 
Since $\mathcal{F}_{k}$ is integrable, we have that $\mathcal{F}_{k}$ is $R$-invariant. We now state a technical lemma, that we shall prove in Section \ref{proofcontd1}. 

\begin{lemma} 
\label{psikub} 
For any $n > k \ge 1$ and for $\left( \mathfrak{B} \times \mu \right)$-almost every $\left( \omega,\, x \right) \in \Sigma_{N}^{+} \times M$, we have 
\[ \phi_{n} \left( \omega,\, x \right)\ \ \le\ \ \sum_{j\, =\, 0}^{n - 1}  \mathcal{F}_{k} \left( \sigma^{j} \omega,\, T_{\pi_{j} \omega} x \right). \]
\end{lemma} 

Making use of $R$-invariance of the function $\mathcal{F}_{k}$ and Lemma \ref{psikub}, we obtain 
\begin{equation} 
\label{intphiuppgeq} 
A\ =\ \lim_{n\, \to\, \infty} \frac{1}{n} \int \phi_{n} \mathrm{d}\left( \mathfrak{B} \times \mu \right) \le\ \int \mathcal{F}_{k} \mathrm{d}\left( \mathfrak{B} \times \mu \right)\ \to\ \int \left( F + \epsilon \right) \mathrm{d}\left( \mathfrak{B} \times \mu \right),\ \ \text{as}\ \ k \to \infty. 
\end{equation} 

Thus, in the case of the sequence of functions $\left\{ \dfrac{\phi_{n}}{n} \right\}_{n\, \ge\, 1}$ being uniformly bounded, we have proved $\displaystyle{\int F \mathrm{d}\left( \mathfrak{B} \times \mu \right) = A}$, refer Equations \eqref{oneway} and \eqref{intphiuppgeq}. 

We now get rid of the uniform bound $- \alpha$. Define for every $\alpha \in \mathbb{R}_{+}$, 
\[ \phi_{n}^{\alpha}\ =\ \max\left\{ \phi_{n}, - n \alpha \right\}\ \ \ \text{and}\ \ \ F^{\alpha}\ =\ \max\left\{ F, - \alpha \right\}. \] 
Then, owing to the subadditivity of the sequence $\left\{ \phi_{n} \right\}_{n\, \ge\, 1}$, we obtain 
\[ \phi_{n + p}^{\alpha} \left( \omega,\, x \right)\ \ \le\ \ \phi_{n}^{\alpha} \left( \omega,\, x \right)\ +\ \phi_{p}^{\alpha } \left( \sigma^{n} \omega,\, T_{\pi_{n} \omega} x \right), \] 
yielding the sequence of functions $\left\{ \phi_{n}^{\alpha} \right\}_{n\, \ge\, 1}$ as a subadditive one with respect to $R$ for any lower bound $- \alpha$. As $\alpha \in \mathbb{R}_{+}$ grows larger, it is easy to note that the sequence $\left\{ \phi_{n}^{\alpha} \left( \omega,\, x \right) \right\}_{\alpha}$ is decreasing for every $ \left( \omega,\, x \right) \in \Sigma_{N}^{+} \times M$. Hence, using the monotone convergence theorem, as can be found in \cite{halmos:1970, rudin:1987}, we get 
\[ \inf_{\alpha\, \in\, \mathbb{R}_{+}} \int \phi_{n}^{\alpha} \mathrm{d}\left( \mathfrak{B} \times \mu \right)\ =\ \int \inf_{\alpha\, \in\, \mathbb{R}_{+}} \phi_{n}^{\alpha} \mathrm{d}\left( \mathfrak{B} \times \mu \right)\ =\ \int \phi_{n} \mathrm{d}\left( \mathfrak{B} \times \mu \right). \] 
Analagously, one can also obtain 
\[ \inf_{\alpha\, \in\, \mathbb{R}_{+}} \int F^{\alpha} \mathrm{d}\left( \mathfrak{B} \times \mu \right)\ =\ \int \inf_{\alpha\, \in\, \mathbb{R}_{+}} F^{\alpha} \mathrm{d}\left( \mathfrak{B} \times \mu \right)\ =\ \int F \mathrm{d}\left( \mathfrak{B} \times \mu \right). \] 
Hence, 
\begin{eqnarray*} 
\inf_{\alpha\, \in\, \mathbb{R}_{+}} \int F^{\alpha} \mathrm{d}\left( \mathfrak{B} \times \mu \right) & = & \inf_{\alpha\, \in\, \mathbb{R}_{+}} \inf_{n\, \ge\, 1} \int \frac{\phi_{n}^{\alpha}}{n} \mathrm{d}\left( \mathfrak{B} \times \mu \right) \\ 
& = & \inf_{n\, \ge\, 1} \inf_{\alpha\, \in\, \mathbb{R}_{+}} \int \frac{\phi_{n}^{\alpha}}{n} \mathrm{d}\left( \mathfrak{B} \times \mu \right) \\ 
& = & \inf_{n\, \ge\, 1} \int \frac{\phi_{n}}{n} \mathrm{d}\left( \mathfrak{B} \times \mu \right). 
\end{eqnarray*} 
Thus, we have proved: 

\begin{lemma} 
\label{integralphil*} 
Let $\left\{ \phi_{n} \right\}$ be a subadditive sequence of measurable functions with respect to $R$ with $\phi_{1}^{+} \in \mathscr{L}^{1} \left( \mathfrak{B} \times \mu \right)$, as considered in Theorem \ref{kingmanomega}. Then, 
\[ \int F \mathrm{d}\left( \mathfrak{B} \times \mu \right)\ \ =\ \ \lim_{n\, \to\, \infty} \int \frac{\phi_{n}}{n} \mathrm{d}\left( \mathfrak{B} \times \mu \right)\ \ =\ \ \inf_{n\, \ge\, 1} \int \frac{\phi_{n}}{n} \mathrm{d}\left( \mathfrak{B} \times \mu \right)\ \ =\ \ A, \] 
where $F = \liminf\limits_{n\, \to\, \infty} \dfrac{\phi_{n}}{n}$. 
\end{lemma} 

We now concentrate on proving $\displaystyle{\int G \mathrm{d}\left( \mathfrak{B} \times \mu \right) \le A}$. As earlier, we again start with the assumption that $\dfrac{\phi_{n}}{n}$ is uniformly bounded below, {\it i.e.}, $\dfrac{\phi_{n}}{n} \left( \omega,\, x \right) \ge - \alpha > - \infty$ for every $\left( \omega,\, x \right) \in \Sigma_{N}^{+} \times M$ for some $\alpha \in \mathbb{R}_{+}$. Fix $k \in \mathbb{Z}_{+}$ and consider for every $n \in \mathbb{Z}_{+}$, a function $\varphi_{n} \left( \omega,\, x \right)$ defined on $\Sigma_{N}^{+} \times M$ and given by 
\[ - \varphi_{n} \left( \omega,\, x \right)\ \ =\ \ \phi_{k} \left( \omega,\, x \right)\ +\ \sum_{j\, =\, 1}^{n - 1} \phi_{k} \left( \sigma^{jk} \omega,\,  T_{\pi_{jk} \omega} x \right). \] 

Since the sequence of functions $\left\{ \varphi_{n} \right\}_{n\, \ge\, 1}$ satisfies the hypothesis of Lemma \ref{integralphil*}, we have, by its assertion that 
\begin{equation} 
\label{apply4.4} 
\int \liminf_{n\, \to\, \infty} \frac{\varphi_{n}}{n} \mathrm{d}\left( \mathfrak{B} \times \mu \right)\ \ =\ \ \lim_{n\, \to\, \infty} \int \frac{\varphi_{n}}{n} \mathrm{d}\left( \mathfrak{B} \times \mu \right). 
\end{equation} 

Owing to our definition of $\varphi_{n}$ and the fact that every function in the sequence $\left\{ \phi_{n} \right\}$ is $R$-invariant, we have 
\begin{equation} 
\label{fomegainv} 
\int \varphi_{n} \mathrm{d}\left( \mathfrak{B} \times \mu \right)\ \ =\ \ - n \int \phi_{k} \mathrm{d}\left( \mathfrak{B} \times \mu \right). 
\end{equation} 
Moreover, using the subadditivity of $\phi_{n}$ with respect to $R$, we have $\varphi_{n} \le - \phi_{nk}$ for all $\left( \omega,\, x \right) \in \Sigma_{N}^{+} \times M$. We now state one more technical lemma that we use to complete the proof of Theorem \ref{kingmanomega}. The proof of this lemma will also be given in the following section. 

\begin{lemma} 
\label{ank=kan} 
Let $\left\{ \phi_{n} \right\}$ be a subadditive sequence of measurable functions with respect to $R$, such that $\phi_{1}^{+} \in \mathscr{L}^{1} (\left( \mathfrak{B} \times \mu \right))$, as considered in Theorem \ref{kingmanomega}. Then, for any fixed $k \in \mathbb{Z}_{+}$, we have 
\[ \limsup_{n\, \to\, \infty} \frac{\phi_{kn}}{n}\ \ =\ \ k \limsup_{n\, \to\, \infty} \frac{\phi_{n}}{n}. \] 
\end{lemma} 

Making use of Lemma \ref{ank=kan}, we now have 
\begin{equation} 
\label{liminfkG} 
\liminf_{n\, \to\, \infty} \frac{\varphi_{n}}{n}\ \le\ - \limsup_{n\, \to\, \infty} \frac{\phi_{nk}}{n}\ =\ - k \limsup_{n\, \to\, \infty} \frac{\phi_{n}}{n}\ =\ -k G. 
\end{equation} 

Integrating the terms in Equation \eqref{liminfkG} with respect to the measure $\left( \mathfrak{B} \times \mu \right)$ and using our findings in Equations in \eqref{apply4.4} and \eqref{fomegainv}, we infer 
\[ - \int \phi_{k} \mathrm{d}\left( \mathfrak{B} \times \mu \right)\ =\ \lim_{n\, \to\, \infty} \int \frac{\varphi_{n}}{n} \mathrm{d}\left( \mathfrak{B} \times \mu \right)\ =\ \int \liminf_{n\, \to\, \infty} \frac{\varphi_{n}}{n} \mathrm{d}\left( \mathfrak{B} \times \mu \right)\ \le\ -k \int G \mathrm{d}\left( \mathfrak{B} \times \mu \right). \] 

Here, making use of Lemma \ref{integralphil*}, after dividing the above inequality by $-k$ and taking limit infimum over $k$, we obtain 
\[ A\ =\ \liminf_{k\, \to\, \infty} \int \frac{\phi_{k}}{k} \mathrm{d}\left( \mathfrak{B} \times \mu \right)\ \ge\ \int G \mathrm{d}\left( \mathfrak{B} \times \mu \right). \] 

We get rid of the uniform bound $\alpha$, as explained earlier. Hence, $F \left( \omega,\, x \right) = G \left( \omega,\, x \right)$ for $\left( \mathfrak{B} \times \mu \right)$-almost every $\left( \omega,\, x \right) \in \Sigma_{N}^{+} \times M$, since $\displaystyle{\int F \mathrm{d}\left( \mathfrak{B} \times \mu \right) = \int G \mathrm{d}\left( \mathfrak{B} \times \mu \right)} = A$. Thus, we define $\phi \equiv F$ on $ \Sigma_{N}^{+} \times M$. Finally, we prove that $\phi$ is $R$-invariant to complete the proof of Theorem \ref{kingmanomega}. 
\begin{eqnarray*} 
\int \left( \phi \circ R \right) \mathrm{d}\left( \mathfrak{B} \times \mu \right) & = & \int \liminf_{n\, \to\, \infty} \frac{\phi_{n} \circ R}{n} \mathrm{d}\left( \mathfrak{B} \times \mu \right)\ \ =\ \ \liminf_{n\, \to\, \infty} \int \frac{\phi_{n} \circ R}{n} \mathrm{d}\left( \mathfrak{B} \times \mu \right) \\ 
& = & \liminf_{n\, \to\, \infty} \int \frac{\phi_{n}}{n} \mathrm{d}\left( \mathfrak{B} \times \mu \right)\ \hspace{+0.85cm} =\ \ \int \liminf_{n\, \to\, \infty} \frac{\phi_{n}}{n} \mathrm{d}\left( \mathfrak{B} \times \mu \right) \\ 
& = & \int \phi \mathrm{d}\left( \mathfrak{B} \times \mu \right). 
\end{eqnarray*} 
\end{proof} 

\section{Proofs of the technical lemmas and Theorem \ref{thm:fkomega}} 
\label{proofcontd1} 

In this section, we will prove Theorem \ref{thm:fkomega}. However, before we go there, as promised in Section \ref{kignmansec}, we provide the proofs of the technical lemmas, namely Lemma \ref{psikub} and Lemma \ref{ank=kan} that we stated and used therein. 

\begin{proof}[of Lemma \ref{psikub}] 
Owing to the subadditivity of the sequence $\phi_{n}$ with respect to $R$, we have 
\begin{eqnarray*} 
\phi_{n} \left( \omega,\, x \right) & \le & \phi_{1} \left( \omega,\, x \right) + \phi_{1} \left( \sigma \omega,\, T_{\pi_{1} \omega} x \right) + \cdots + \phi_{1} \left( \sigma^{n - 1} \omega,\, T_{\pi_{n - 1} \omega} x \right) \\ 
& = & \phi_{1} \left( \omega,\, x \right)\ +\ \sum_{1\, \le\, j\, \le\, n - 1\; :\; \left( \sigma^{j} \omega,\, T_{\pi_{j} \omega} x \right)\, \in\, \mathcal{H}_{k}} \phi_{1} \left( \sigma^{j} \omega,\, T_{\pi_{j} \omega} x \right) \\ 
& & \hspace{+1.8cm} + \sum_{1\, \le\, j\, \le\, n - 1\; :\; \left( \sigma^{j} \omega,\, T_{\pi_{j} \omega} x \right)\, \notin\, \mathcal{H}_{k}} \phi_{1} \left( \sigma^{j} \omega,\, T_{\pi_{j} \omega} x \right). 
\end{eqnarray*} 
Using the definitions of the set $\mathcal{H}_{k}$ and the function $\mathcal{F}_{k}$, as written in Equations \eqref{hkw} and \eqref{fkw}, we obtain 
\begin{eqnarray*} 
\phi_{n} \left( \omega,\, x \right) & \le & \phi_{1} \left( \omega,\, x \right) + \sum_{1\, \le\, j\, \le\, n - 1\; :\; \left( \sigma^{j} \omega,\, T_{\pi_{j} \omega} x \right)\, \in\, \mathcal{H}_{k}} \left( F \left( \sigma^{j} \omega,\, T_{\pi_{j} \omega} x \right) + \epsilon \right) \\ 
& & \hspace{+1.8cm} + \sum_{1\, \le\, j\, \le\, n - 1\; :\; \left( \sigma^{j} \omega,\, T_{\pi_{j} \omega} x \right)\, \notin\, \mathcal{H}_{k}} \phi_{1} \left(\sigma^{j} \omega,\, T_{\pi_{j} \omega} x \right) \\ 
& \le & \mathcal{F}_{k} \left( \omega,\, x \right)\ + \ \sum_{1\, \le\, j\, \le\, n - 1\; :\; \left( \sigma^{j} \omega,\, T_{\pi_{j} \omega} x \right)\, \in\, \mathcal{H}_{k}} \left( F \left( \sigma^{j} \omega,\, T_{\pi_{j} \omega} x \right) + \epsilon \right) \\ 
& & \hspace{+1.8cm} + \sum_{1\, \le\, j\, \le\, n - 1\; :\; \left( \sigma^{j} \omega,\, T_{\pi_{j} \omega} x \right)\, \notin \, \mathcal{H}_{k}} \mathcal{F}_{k} \left( \sigma^{j} \omega,\, T_{\pi_{j} \omega} x \right) \\ 
& = & \mathcal{F}_{k} \left( \omega,\, x \right) + \sum_{j\, =\, 1}^{n - 1} \mathcal{F}_{k} \left( \sigma^{j} \omega,\, T_{\pi_{j} \omega} x\right). 
\end{eqnarray*} 
\end{proof} 

\begin{proof}[of Lemma \ref{ank=kan}] 

Fix $k \in \mathbb{Z}_{+}$ and note that $\left\{ \phi_{kn} \right\}_{n\, \ge\, 1}$ is a subsequence of $\left\{ \phi_{n} \right\}_{n\, \ge\, 1}$ and thus, 
\[ \limsup_{n\, \to\, \infty} \frac{\phi_{kn}}{kn}\ \ \le\ \ \limsup_{n\, \to\, \infty} \frac{\phi_{n}}{n}. \] 

To prove the other way around, let $n = k q_{n} + r_{n}$ where $q_{n} \in \mathbb{Z}_{+}$ and $r_{n} \in \left\{ 0, 1, 2, \cdots, k - 1 \right\}$. Then, 
\[ \phi_{n} \left( \omega,\, x \right)\ \le\ \phi_{k q_{n}} \left( \omega,\, x \right) + \phi_{r_{n}} \left( \sigma^{k q_{n}} \omega,\, T_{\pi_{k q_{n}} \omega} x \right). \] 
Define a function $\psi : \Sigma_{N}^{+} \times  M \longrightarrow \mathbb{R}$ by $\psi \left( \sigma^{j_{0} q_{n}} \omega,\, T_{\pi_{j_{0} q_{n}} \omega} x \right) = \max\limits_{1\, \le\, j\, \le\, k - 1} \left\{ \phi_{j}^{+} \left( \sigma^{j q_{n}} \omega,\, T_{\pi_{j q_{n}} \omega} x \right) \right\}$. Then, 
\[ \phi_{n} \left( \omega,\, x \right)\ \le\ \phi_{k q_{n}} \left( \omega,\, x \right) + \psi \left( \sigma^{j_{0} q_{n}} \omega,\, T_{\pi_{j_{0} q_{n}} \omega} x \right). \] 

Taking limit supremum after dividing the above inequality by $n$, we obtain 
\[ \limsup_{n\, \to\, \infty} \frac{1}{n} \phi_{n} \left( \omega,\, x \right)\ \le\ \limsup_{n\, \to\, \infty} \frac{1}{n} \phi_{k q_{n}} \left( \omega,\, x \right) + \limsup_{n\, \to\, \infty} \frac{1}{n} \psi \left( \sigma^{j_{0} q_{n}} \omega,\, T_{\pi_{j_{0} q_{n}} \omega} x \right). \] 

We now prove that the second term in the right hand side of the above inequality vanishes for $\left( \mathfrak{B} \times \mu \right)$-almost every $\left( \omega,\, x \right) \in \Sigma_{N}^{+} \times M$, by using the Borel-Cantelli lemma, that we now state. 

\begin{lemma}[Borel-Cantelli Lemma]\cite{parthasarathy:1977} 
\label{borelcanteli}
For a sequence $\left\{ B_{n} \right\}$ in the $\sigma$-algebra of a probability space $\left( M, \mu \right)$ that satisfies $\sum\limits_{n\, \ge\, 1} \mu \left( B_{n} \right) < \infty$, we have $\mu \left( \limsup\limits_{n\, \to\, \infty} B_{n} \right) = 0$. 
\end{lemma} 

Note that by definition $\psi \in \mathscr{L}^{1} \left( \mathfrak{B} \times \mu \right) $ and hence, $\psi$ is $R$-invariant. Fix $\epsilon > 0$, define a sequence of sets $B_{n} = \big\{ \left( \omega,\, x \right) \in \Sigma_{N}^{+} \times M : \left| \psi \left( \sigma^{n} \omega,\, T_{\pi_{n} \omega} x \right) \right| \ge n \epsilon \big\}$ and consider the series 
\begin{eqnarray*} 
\sum_{n\, \ge\, 1} \mu \left( B_{n} \right) & = & \sum_{n\, \ge\, 1} \mu \left( \left\{ \left( \omega,\, x \right) \in \Sigma_{N}^{+} \times M\ :\ \left| \psi \left( \omega,\, x \right) \right| \ge n \epsilon \right\} \right) \\ 
& = & \sum_{n\, \ge\, 1} \sum_{k\, \ge\, n} \mu \left( \left\{ \left( \omega,\, x \right) \in \Sigma_{N}^{+} \times M\ :\ k \le \frac{\left| \psi \left( \omega,\, x \right) \right|}{\epsilon} < k + 1 \right\} \right) \\ 
& = & \sum_{k\, \ge\, 1} k \mu \left( \left\{ \left( \omega,\, x \right) \in \Sigma_{N}^{+} \times M\ :\ k \le \frac{\left| \psi \left( \omega,\, x \right) \right|}{\epsilon} < k + 1 \right\} \right) \\ 
& \le & \frac{1}{\epsilon} \int \left| \psi \right| \mathrm{d}\left( \mathfrak{B} \times \mu \right) \\ 
& < & \infty. 
\end{eqnarray*}

Thus, appealing to the Borel-Cantelli lemma \ref{borelcanteli}, we obtain $\mu (B) = 0$, where 
\[ B\ \ =\ \ \limsup\limits_{n\, \to\, \infty} B_{n}\ \ =\ \ \Big\{ \left( \omega,\, x \right) \in \Sigma_{N}^{+} \times M\ :\left| \psi \left( \sigma^{n} \omega,\, T_{\pi_{n} \omega} x \right) \right| \ge n \epsilon\ \text{for infinitely many}\ n \in \mathbb{Z}_{+} \Big\}. \] 

Thus, if $\left( \tau,\, y \right) \notin B$, then there exists a threshold, say $n_{0}$ such that $\left| \psi \left( \sigma^{n} \tau,\, T_{\pi_{n} \tau} y \right) \right| < n \epsilon$ for all $n > n_{0}$. Choosing $\epsilon$ arbitrarily small, we then obtain $\dfrac{\left| \psi \left( \sigma^{n} \tau,\, T_{\pi_{n} \tau} y \right) \right|}{n}$ and therefore, $\dfrac{\psi \circ R^{n}}{n}$, to be sufficiently close to $0$ on $\left( \Sigma_{N}^{+} \times M \right) \setminus B$. Hence, 
\[ \limsup_{n\, \to\, \infty} \frac{1}{n} \phi_{n} \left( \omega,\, x \right)\ \le\ \limsup_{n\, \to\, \infty} \frac{1}{k q_{n}} \phi_{k q_{n}} \left( \omega,\, x \right)\ \le\ \frac{1}{k} \limsup_{n\, \to\, \infty} \frac{1}{n} \phi_{nk} \left( \omega,\, x \right). \] 
\end{proof} 

\begin{proof}[of Theorem \ref{thm:fkomega}] 

Define a sequence of functions $\left\{ \phi_{n} \right\}_{n\, \ge\, 1};\; \phi_{n} : \Sigma_{N}^{+} \times M \longrightarrow [- \infty, \infty)$ by 
\[ \phi_{n} \left( \omega,\, x \right)\ \ =\ \ \log \left\| \prod_{j\, =\, 0}^{n - 1}  L \left( T_{\pi_{j} \omega} x \right) \right\|_{{\rm op}}. \] 
Then, an application of Theorem \ref{kingmanomega} completes the proof for the existence, $R$-invariance and the value of the integral as mentioned in the statement of the theorem for $\lambda_{+}$. 

Similarly, considering the sequence of functions $\left\{ \psi_{n} \right\}_{n\, \ge\, 1};\; \psi_{n} : \Sigma_{N}^{+} \times M \longrightarrow [- \infty, \infty)$ defined by 
\[ \psi_{n} \left( \omega,\, x \right)\ \ =\ \ \log \left\| \left( \prod_{j\, =\, 0}^{n - 1} L \left( T_{\pi_{j} \omega} x \right) \right)^{-1} \right\|_{{\rm op}}^{-1}, \] 
completes the proof for $\lambda_{-}$. 
\end{proof} 

\section{Proofs of Theorem \ref{kingmanomegafree} and Theorem \ref{thm:fkomegafree}} 
\label{proofcontd2} 

In this section, we prove Theorem \ref{kingmanomegafree} and Theorem \ref{thm:fkomegafree}. Even though the proofs may appear to follow exactly along the lines of the proofs of Theorem \ref{kingmanomega} and Theorem \ref{thm:fkomega} respectively, one needs to remember that the map $\mathfrak{S}$ that determines simultaneous dynamics is a set-valued map defined on subsets of $M \times \mathbb{R}^{d}$. Further, the dynamics due to $\mathfrak{S}$ as well as $S$ accommodates the totality of all branches of all points in the considered set. One obvious account that allows us to immediately see this change is the presence of a factor of $n N^{n}$ in the definition of the Lyapunov exponents, as written in Equation \eqref{lambdaplusomegafree} in the case of simultaneous dynamics, in place of the polynomial factor $n$ in the definition of the Lyapunov exponents, as written in Equation \eqref{lambdaplmiomega} in the case of random dynamics. 

\begin{proof}[of Theorem \ref{kingmanomegafree}]  
As earlier, we begin with the hypothesis that states that $\left\{ \Phi_{n} \right\}_{n\, \ge\, 1}$ is a subadditive sequence of measurable functions with respect to $S$. This yields 
\[ \Phi_{n} (x)\ \ \le\ \ \Phi_{1} (x) + \left( \Phi_{1} \circ S \right) (x) + \left( \Phi_{1} \circ S^{2} \right) (x) + \cdots + \left( \Phi_{1} \circ S^{n - 1} \right) (x), \] 
where $\left( \Phi_{1} \circ S^{k} \right)$ is as defined in Equation \eqref{composition}. Further, 
\begin{eqnarray*} 
\int \Phi_{n} \mathrm{d}\mu & \le & \int \Phi_{1} \mathrm{d}\mu + \sum_{\left\{ \pi_{1} \omega\; :\; \omega\, \in\, \Sigma_{N}^{+} \right\}} \int \left( \Phi_{1} \circ T_{\pi_{1} \omega} \right) \mathrm{d}\mu  + \sum_{\left\{ \pi_{2} \omega\; :\; \omega\, \in\, \Sigma_{N}^{+} \right\}} \int \left( \Phi_{1} \circ T_{\pi_{2} \omega} \right) \mathrm{d}\mu \\ 
& & \hspace{+4cm} + \cdots + \sum_{\left\{ \pi_{n - 1} \omega\; :\; \omega\, \in\, \Sigma_{N}^{+} \right\}} \int \left( \Phi_{1} \circ T_{\pi_{n - 1} \omega} \right) \mathrm{d}\mu \\ 
& \le & \left( 1 + N + N^{2} + \cdots + N^{n - 1} \right) \int \Phi_{1} \mathrm{d}\mu \\ 
& = & \frac{N^{n} - 1}{N - 1} \int \Phi_{1} \mathrm{d}\mu. 
\end{eqnarray*} 

Hence, defining $\displaystyle{\mathfrak{a}_{n} = \frac{1}{N^{n}} \int \Phi_{n} \mathrm{d}\mu}$, we note that $\left\{ \mathfrak{a}_{n} \right\}$ is a subadditive sequence and thus, by Lemma \ref{lim=inf}, we have $\displaystyle{\mathfrak{a} = \lim_{n\, \to\, \infty} \int \frac{\Phi_{n}}{n N^{n}} \mathrm{d}\mu = \inf_{n\, \ge\, 1} \int \frac{\Phi_{n}}{n N^{n}} \mathrm{d}\mu}$. We now define 
\[ \mathpzc{F}(x)\ \ =\ \ \liminf_{n\, \to\, \infty} \frac{\Phi_{n} (x)}{n N^{n}}\ \ \ \text{and}\ \ \ \mathpzc{G}(x)\ \ =\ \ \limsup_{n\, \to\, \infty} \frac{\Phi_{n} (x)}{n N^{n}}. \] 
Proving $\displaystyle{\int \mathpzc{G} \mathrm{d}\mu \le \mathfrak{a} \le \int \mathpzc{F} \mathrm{d}\mu}$ completes the proof of Theorem \ref{kingmanomegafree}. Our next step is to consider $\dfrac{\Phi_{n} (x)}{n N^{n}} \ge - \beta > - \infty$ for every $x \in M$ for some $\beta \in \mathbb{R}_{+}$. Then, applying Fatou's lemma, as may be found in \cite{halmos:1970, rudin:1987}, we get $\displaystyle{\int \mathpzc{F} \mathrm{d}\mu \le \mathfrak{a}}$. Further, fixing $\epsilon > 0$, we define a nested exhausting sequence of subsets of $M$, say $\mathpzc{H}_{k}$ given by 
\[ \mathpzc{H}_{k}\ \ =\ \ \left\{ x \in M\ :\ \Phi_{j} (x) \le j \left( \mathpzc{F}(x) + \epsilon \right)\ \text{for some}\ 1 \le j \le k \right\}. \] 
Then, defining $\mathscr{F}_{k} : M \longrightarrow \mathbb{R}$ as 
\[ \mathscr{F}_{k} (x)\ \ =\ \ \begin{cases} \mathpzc{F}(x) + \epsilon & \text{if}\ x \in \mathpzc{H}_{k} \\ \Phi_{1}(x) & \text{if}\ x \in M \setminus \mathpzc{H}_{k}, \end{cases} \]
and using the monotone convergence theorem,
we obtain $\displaystyle{\int \mathscr{F}_{k} \mathrm{d}\mu \to \int \mathpzc{F} \mathrm{d}\mu + \epsilon}$. 

We now state a lemma, analogous to Lemma \ref{psikub} in the current context. 

\begin{lemma} 
\label{phikupomegafree} 
For any $n > k \ge 1$ and for $\mu$-almost every $x \in M$, we have 
\[ \Phi_{n} (x)\ \ \le\ \ \mathscr{F}_{k} (x) + \sum_{j\, =\, 1}^{n - 1} \sum_{\left\{ \pi_{j} \omega\; :\; \omega\, \in \Sigma_{N}^{+} \right\}} \left( \mathscr{F}_{k} \circ T_{\pi_{j} \omega} \right) (x). \]
\end{lemma} 

Economising on the length of the paper, we do not write the proof of Lemma \ref{phikupomegafree}, and merely remark that the arguments in this proof can be written {\it mutatis mutandis} to the proof of Lemma \ref{psikub}. Now making use of Lemma \ref{phikupomegafree}, we get 
\begin{eqnarray*} 
\mathfrak{a} & \le & (N - 1) \lim_{n\, \to\, \infty} \frac{1}{n N^{n}} \int \Phi_{n} \mathrm{d}\mu \\ 
& \le & (N - 1) \lim_{n\, \to\, \infty} \frac{\left( 1 + N + N^{2} + \cdots + N^{n - 1} \right)}{n N^{n}} \int \mathscr{F}_{k} \mathrm{d}\mu \\ 
& \le & \lim_{n\, \to\, \infty} \frac{N^{n} - 1}{N^{n}} \int \mathscr{F}_{k} \mathrm{d}\mu \\ 
& \to & \int \mathpzc{F} \mathrm{d}\mu + \epsilon,\ \ \ \text{as}\ k \to \infty. 
\end{eqnarray*} 

To get rid of the uniform bound $- \beta$, we define for $\beta \in \mathbb{R}_{+}$, 
\[ \Phi_{n}^{\beta}\ \ =\ \ \max\left\{ \Phi_{n}, - n N^{n} \beta \right\}\ \ \ \ \text{and}\ \ \ \ \mathpzc{F}^{\beta}\ \ =\ \ \max\left\{ \mathpzc{F}, - \beta \right\}. \] 
Again going through the lines of the proof of Theorem \ref{kingmanomega}, we arrive at 
\begin{equation} 
\label{intphin=a} 
\mathfrak{a}\ \ =\ \ \lim_{n\, \to\, \infty} \int \frac{\Phi_{n}}{n N^{n}} \mathrm{d}\mu\ \ =\ \ \inf_{n\, \ge\, 1} \int \frac{\Phi_{n}}{n N^{n}} \mathrm{d}\mu. 
\end{equation} 

We finally prove $\displaystyle{\int \mathpzc{G} \mathrm{d}\mu \le \mathfrak{a}}$. As in the proof of Theorem \ref{kingmanomega}, we start with the assumption $\dfrac{\Phi_{n}(x)}{n N^{n}} \ge - \beta > - \infty$ for every $x \in M$ and for some fixed $\beta \in \mathbb{R}_{+}$. Here, for a fixed $k \in \mathbb{Z}_{+}$ and for every $n \in \mathbb{Z}_{+}$, we define $\xi_{n}$ on $M$ as 
\begin{equation}
\label{xinomegafree}    
- \frac{1}{n} \xi_{n} (x)\ =\ \Phi_{k} (x)\ +\ \sum_{j\, =\, 1}^{n - 1} \sum_{\left\{ \pi_{jk} \omega\; :\; \omega\, \in\, \Sigma_{N}^{+} \right\}} \Phi_{k} \left( T_{\pi_{jk} \omega} x \right). 
\end{equation}

Thus, observing that $\xi_{n}$ is a subadditive sequence of measurable functions on $M$ with respect to $S$ and making use of the idea that gave Equation \eqref{intphin=a}, we now obtain 
\begin{equation} 
\label{firstpartomegafree} 
\int \liminf_{n\, \to\, \infty} \frac{\xi_{n}}{n N^{n}} \mathrm{d}\mu\ \ =\ \ \lim_{n\, \to\, \infty} \int \frac{\xi_{n}}{n N^{n}} \mathrm{d}\mu. 
\end{equation} 
Also, $\dfrac{\xi_{n}}{n} \le - \Phi_{nk}$ for every $x \in M$ and thus, $\dfrac{\xi_{n}}{n N^{n}} \le - \dfrac{\Phi_{nk}}{N^{n}} \le - \dfrac{\Phi_{nk}}{n N^{n}}$ for every $x \in M$. 

We now state a lemma analogous to Lemma \ref{ank=kan} and provide a short proof of the same. 

\begin{lemma} 
\label{akn=kanomegafree} 
Let $\left\{ \Phi_{n} \right\}$ be a subadditive sequence of measurable functions with respect to $S$, such that $\Phi_{1}^{+} \in \mathscr{L}^{1} (\left( \mu \right))$, as considered in Theorem \ref{kingmanomegafree}. Then, for any fixed $k \in \mathbb{Z}_{+}$, we have 
\[ \limsup_{n\, \to\, \infty} \frac{\Phi_{kn}}{n N^{n}}\ \ =\ \ k N^{k} \limsup_{n\, \to\, \infty} \frac{\Phi_{n}}{n N^{n}}. \] 
\end{lemma} 

\begin{proof} 
Observe that $\left\{ \dfrac{\Phi_{nk}}{nk N^{nk}} \right\}$ is a subsequence of the sequence $\left\{ \dfrac{\Phi_{n}}{n N^{n}} \right\}$. This implies 
\[ \limsup_{n\, \to\, \infty} \frac{\Phi_{nk}}{nk N^{nk}}\ \ \le\ \ \limsup_{n\, \to\, \infty} \frac{\Phi_{n}}{n N^{n}}. \] 
Moreover, for large $n$, we have $nk N^{nk} \ge nk N^{n + k} \ge nk N^{n} \ge n N^{n}$. Thus, 
\[ \limsup_{n\, \to\, \infty} \frac{\Phi_{nk}}{nk N^{nk}}\ \ \le\ \ \limsup_{n\, \to\, \infty} \frac{\Phi_{nk}}{nk (N^{n + k})}\ \ \le\ \ \limsup_{n\, \to\, \infty} \frac{\Phi_{n}}{n N^{n}}, \] 
where the last inequality is obtained by using the subadditivity of the sequence $\Phi_{n}$. Thus, 
\[ \limsup_{n\, \to\, \infty} \frac{\Phi_{nk}}{n N^{n}}\ \ \le\ \ k N^{k} \limsup_{n\, \to\, \infty} \frac{\Phi_{n}}{n N^{n}}, \] 
proving one way inequality. 

On the other hand, consider $n = k q_{n} + r_{n}$ where $r_{n} \in \{ 0, 1, \cdots, k - 1 \}$. Then, by employing the definition of subadditivity, we have $\Phi_{n} \le \Phi_{k q_{n}} + \sum\limits_{\left\{ \pi_{k q_{n}} \omega\; :\; \omega\, \in\, \Sigma_{N}^{+} \right\}} \Phi_{r_{n}} \circ T_{\pi_{k q_{n}} \omega}$, where we take $\Phi_{0} \equiv 0$. Then, 
\begin{eqnarray*} 
\frac{\Phi_{n}}{n N^{n}} & \le & \frac{\Phi_{k q_{n}}}{n N^{n}} + \frac{1}{n N^{n}} \sum_{\left\{ \pi_{k q_{n}} \omega\; :\; \omega\, \in\, \Sigma_{N}^{+} \right\}} \Phi_{r_{n}} \circ T_{\pi_{k q_{n}} \omega} \\ 
& = & \frac{1}{k q_{n} + r_{n}} \frac{\Phi_{k q_{n}}}{N^{k q_{n} + r_{n}}} + \frac{1}{n N^{n}} \sum_{\left\{ \pi_{k q_{n}} \omega\; :\; \omega\, \in\, \Sigma_{N}^{+} \right\}} \Phi_{r_{n}} \circ T_{\pi_{k q_{n}} \omega} \\ 
& \le & \frac{\Phi_{k q_{n}}}{k q_{n} N^{k q_{n}}} + \frac{1}{n N^{n}} \sum_{\left\{ \pi_{k q_{n}} \omega\; :\; \omega\, \in\, \Sigma_{N}^{+} \right\}} \Phi_{r_{n}} \circ T_{\pi_{k q_{n}} \omega} \\ 
& \le & \frac{\Phi_{k q_{n}}}{k q_{n} N^{k + q_{n}}} + \frac{1}{k q_{n} N^{k q_{n}}} \sum_{\left\{ \pi_{k q_{n}} \omega\; :\; \omega\, \in\, \Sigma_{N}^{+} \right\}} \Phi_{r_{n}} \circ T_{\pi_{k q_{n}} \omega}\ \ \text{for large}\ n, 
\end{eqnarray*} 
since $\left( k q_{n} + r_{n} \right) N^{k q_{n} + r_{n}} \ge k q_{n} N^{k q_{n}}$. To complete the proof of the lemma, we will now assume that the second term goes to $0$ as $n \to \infty$ for $\mu$-almost every $x$. We will prove the same after finishing the following step. Note that 
\[ \limsup_{n\, \to\, \infty} \frac{\Phi_{n}}{n N^{n}}\ \ \le\ \ \frac{1}{k N^{k}} \limsup_{q_{n}\, \to\, \infty} \frac{\Phi_{k q_{n}}}{q_{n} N^{q_{n}}}\ \ =\ \ \frac{1}{k N^{k}} \limsup_{n\, \to\, \infty} \frac{\Phi_{k n}}{n N^{n}}, \] 
proving the other way inequality, and thus the lemma. 

Now, we prove the intermediary step, where we have assumed 
\[ \lim_{n\, \to\, \infty} \frac{1}{k q_{n} N^{k q_{n}}} \sum_{\left\{ \pi_{k q_{n}} \omega\; :\; \omega\, \in\, \Sigma_{N}^{+} \right\}} \Phi_{r_{n}} \circ T_{\pi_{k q_{n}} \omega} (x) = 0\ \ \text{for}\ \mu\text{-almost every}\ x. \] 

Define $\Psi = \max \left\{ \Phi_{0}, \Phi_{1}, \cdots, \Phi_{k - 1} \right\}$, where $k \in \mathbb{Z}_{+}$ is as fixed in the statement of Lemma \ref{akn=kanomegafree}. Since $\Phi_{1}^{+} \in \mathscr{L}^{1} (\mu)$, we obtain by subadditivity of the sequence that $\Phi_{n}$ is integrable and thus, $\Psi$ is integrable too. Fix $\epsilon > 0$ and for $m = k q_{n}$ define 
\[ \mathcal{B}_{m}\ \ =\ \ \left\{ x \in M\ :\ \left| \sum_{\left\{ \pi_{m} \omega\; :\; \omega\, \in\, \Sigma_{N}^{+} \right\}} \Psi \circ T_{\pi_{m} \omega} (x) \right| \ge m N^{m} \epsilon \right\}. \] 
Then, 
\[ \sum_{m\, \ge\, 1} \mu \left( \mathcal{B}_{m} \right)\ =\ \sum_{m\, \ge\, 1} \mu \left( \left\{ x \in M\ :\ \left| \Psi (x) \right| \ge m \epsilon \right\} \right)\ \le\ \frac{1}{\epsilon} \int | \Psi | \mathrm{d}\mu\ <\ \infty, \]
as in the proof of Lemma \ref{ank=kan}. Using the Borel-Cantelli lemma, we conclude that $\mu \left( \mathcal{B} \right) = 0$, where 
\begin{eqnarray*} 
\mathcal{B} & = & \limsup\limits_{m\, \to\, \infty} \mathcal{B}_{m} \\ 
& = & \left\{ x \in M\ :\ \left| \sum_{\left\{ \pi_{m} \omega\; :\; \omega\, \in\, \Sigma_{N}^{+} \right\}} \Psi \circ T_{\pi_{m} \omega} (x) \right| \ge m N^{m} \epsilon\ \text{for infinitely many}\ m \in \mathbb{Z}_{+} \right\}. 
\end{eqnarray*} 

Thus, if $y \notin \mathcal{B}$, then there exists a threshold, say $m_{0}$ such that 
\[ \left| \sum_{\left\{ \pi_{m} \omega\; :\; \omega\, \in\, \Sigma_{N}^{+} \right\}} \Psi \circ T_{\pi_{m} \omega} (y) \right| < m N^{m} \epsilon\ \ \text{for all}\ \ m > m_{0}. \] 
Finally, we choose $\epsilon$ arbitrarily small to conclude the proof of this lemma, along the lines of the proof of Lemma \ref{ank=kan}. 
\end{proof} 

Making use of Lemma \ref{akn=kanomegafree}, we now obtain 
\begin{equation} 
\label{bddforg} 
\liminf_{n\, \to\, \infty} \frac{\xi_{n}}{n N^{n}}\ \le\ - \limsup_{n\, \to\, \infty} \frac{\Phi_{nk}}{n N^{n}}\ =\ - k N^{k} \limsup_{n\, \to\, \infty} \frac{\Phi_{n}}{n N^{n}}\ =\ - k N^{k} \mathpzc{G}. 
\end{equation} 

Integrating the function $\xi_{n}$ as defined in Equation \eqref{xinomegafree}, with respect to the measure $\mu$, we obtain 
\begin{equation} 
\label{hjgu} 
- \int \Phi_{k} \mathrm{d}\mu\ \ \le\ \ \frac{N^{k} - 1}{n N^{nk}} \int \xi_{n} \mathrm{d}\mu\ \ \le\ \ \frac{1}{n N^{n}} \int \xi_{n} \mathrm{d}\mu. 
\end{equation} 
Employing Equation \eqref{firstpartomegafree} and Inequality \eqref{bddforg} in Inequality \eqref{hjgu} and applying limits, we get 
\[ - \int \Phi_{k} \mathrm{d}\mu\ \le\ \lim_{n\, \to\, \infty} \int \frac{\xi_{n}}{n N^{n}} \mathrm{d}\mu\ =\ \int \liminf_{n\, \to\, \infty} \frac{\xi_{n}}{n N^{n}} \mathrm{d}\mu\ \le\ - k N^{k} \int \mathpzc{G} \mathrm{d}\mu. \] 
Dividing the above inequality by $k N^{k}$ and taking limit infimum as $k \to \infty$, we get 
\[ \mathfrak{a}\ \ =\ \ \liminf_{k\, \to\, \infty} \int \frac{\Phi_{k}}{k N^{k}} \mathrm{d}\mu\ \ \ge\ \ \int \mathpzc{G} \mathrm{d}\mu. \] 
One can remove the uniform bound $\beta$ by following the same lines as explained in the proof of Theorem \ref{kingmanomega}. Further, we observe that 
\begin{eqnarray*} 
\int \left( \Phi \circ S \right) \mathrm{d}\mu & = & \int \liminf_{n\, \to\, \infty} \frac{\Phi_{n} \circ S}{n N^{n}} \mathrm{d}\mu\ \ =\ \ \liminf_{n\, \to\, \infty} \int \frac{\Phi_{n} \circ S}{n N^{n}} \mathrm{d}\mu \\ 
& = & \liminf_{n\, \to\, \infty} \int \frac{N \Phi_{n}}{n N^{n}} \mathrm{d}\mu\ \hspace{+0.4cm} =\ \ \int \liminf_{n\, \to\, \infty} \frac{N \Phi_{n}}{n N^{n}} \mathrm{d}\mu\ \hspace{+0.5cm} =\ \ \int N \Phi \mathrm{d}\mu, 
\end{eqnarray*} 
thus, proving the $S$-invariance of the limit function $\Phi$ and thereby, completing the proof. 
\end{proof} 

We finally conclude the paper with the proof of Theorem \ref{thm:fkomegafree}. 

\begin{proof}[of Theorem \ref{thm:fkomegafree}] 
Define two sequences of functions $\Phi_{n}$ and $\Psi_{n}$ from $M$ to $[- \infty, \infty)$ by 
\begin{eqnarray*} 
\Phi_{n} (x) & = & \sum_{\left\{ \pi_{n} \omega\; :\; \omega\, \in\, \Sigma_{N}^{+} \right\}} \log \left\| \prod_{j\, =\, 0}^{n - 1} L \left( T_{\pi_{j} \omega} (x) \right) \right\|_{{\rm op}}\ \ \ \ \text{and} \\ 
\Psi_{n} (x) & = & \sum_{\left\{ \pi_{n} \omega\; :\; \omega\, \in\, \Sigma_{N}^{+} \right\}} \log \left\| \left( \prod_{j\, =\, 0}^{n - 1} L \left( T_{\pi_{j} \omega} (x) \right) \right)^{-1} \right\|_{{\rm op}}^{-1}. 
\end{eqnarray*} 
Then, the proof of Theorem \ref{thm:fkomegafree} is complete, when viewed as a corollary to Theorem \ref{kingmanomegafree}. 
\end{proof} 
\bigskip 

{\bf Acknowledgements:} The authors are thankful to Pablo G. Barrientos and Dominique Malicet for identifying an error in the proof of Theorem \ref{kingmanomegafree}, in an earlier version of the paper. 
\bigskip 

{\bf Data Availability statement:} The authors declare that the data supporting the findings of this study are available within the paper. 
 
\bigskip 
\bigskip 

\bigskip 

\bigskip 



\bigskip 


\begin{thebibliography}{99} 

\bibitem{as:1978} 
{\sc Akcoglu, M. A.} and {\sc Sucheston, L.}, ``A ratio ergodic theorem for superadditive processes", \emph{Z. Wahrsch. Verw. Gebiete}, {\bf 44} (1978) 269 - 278. 

\bibitem{al:1998}
{\sc Arnold, L.}, ``Random dynamical systems", \emph{Springer Monographs in Mathematics}, Springer, Berlin, (1998). 

\bibitem{aev:pp2023} 
{\sc Avila, A., Eskin, A.} and {\sc Viana, M.}, ``Continuity of the Lyapunov exponents of random matrix products", \texttt{arXiV:2305.06009} (2023). 

\bibitem{bm:2020} {\sc Barrientos, P. G.} and {\sc Malicet, D.}, ``Extremal exponents of random products of conservative diffeomorphisms", \emph{Math. Z.}, {\bf 296} (2020) 1185 - 1207. 

\bibitem{bm:pp} {\sc Barrientos, P. G.} and {\sc Malicet, D.}, ``Mostly contracting random maps", \texttt{arXiV:2412:03729} (2024). 


\bibitem{bv:2017} 
{\sc Bocker-Neto, C.} and {\sc Viana, M.}, ``Continuity of Lyapunov exponents for random two-dimensional matrices", \emph{Ergodic Theory Dynam. Systems}, {\bf 37} (2017) 1413 - 1442. 



 





\bibitem{coor:2015}
{\sc Coornaert, M.}, ``Topological dimension and dynamical systems", \emph{Universitext}, (2015).


 





\bibitem{dk:2016} 
{\sc Duarte, P. M. N. R. D.} and {\sc Klein, S.}, ``Lyapunov exponents of linear cocycles", \emph{Atlantis Studies in Dynamical Systems}, {\bf 3}, Atlantis Press, Paris, 2016. 


\bibitem{furs:1960}
{\sc Furstenberg, H.} and {\sc Kesten, H.}, ``Products of random matrices", \emph{Ann. Math. Statist.}, {\bf 31} (1960) 457 - 469.

\bibitem{aswin:2024}
{\sc Gopakumar, A., Rajasekar, K.} and {\sc Sridharan, S.}, ``Simultaneous action of finitely many interval maps: some dynamical and statistical properties", \emph{Real Anal. Exchange}, {\bf 49} (2024) 13 - 66. 


\bibitem{halmos:1970}
{\sc Halmos, P. R.}, ``Measure Theory", \emph{D. Van Nostrand Co.}, Inc., New York, (1950).


\bibitem{kingman:1968}
{\sc Kingman, J. F. C.}, ``The ergodic theorem of subadditive stochastic processes", \emph{J. Royal Statist. Soc.} Ser. B, {\bf 30} (1968) 499 - 510. 


\bibitem{lq:1995} 
{\sc Liu P.-D.} and {\sc Qian, M.}, ``Smooth ergodic theory of random dynamical systems", \emph{Lecture Notes in Mathematics}, (1995). 


\bibitem{parthasarathy:1977}
{\sc Parthasarathy, K. R.}, ``Introduction to probability and measure", \emph{Texts Read. Math.}, {\bf 33}, Hindustan Book Agency, (2005). 




\bibitem{poli:1998}
{\sc Pollicott, M.} and {\sc Yuri, M.}, ``Dynamical systems and ergodic theory", \emph{London Math. Soc. Stud. Texts}, {\bf 40} (1998).

\bibitem{rudin:1987}
{\sc Rudin, W.}, ``Real and complex Analysis", \emph{McGraw-Hill Book Co., New York}, (1987).



\bibitem{tv:2020} 
{\sc Tall, E. H. Y.} and {\sc Viana, M.}, ``Moduli of continuity for the Lyapunov exponents of random ${\rm GL}(2)$-cocycles", \emph{Trans. Amer. Math. Soc.}, {\bf 373} (2020) 1343 - 1383. 


\bibitem{mv:2014}
{\sc Viana, M.}, ``Lectures on Lyapunov exponents", \emph{Cambridge Stud. Adv. Math.}, {\bf 145} (2014).

\bibitem{pw:2000}
{\sc Walters, P.},  ``An introduction to ergodic theory", \emph{Grad. Texts in Math.}, (1982). 


\end{thebibliography}
\end{document}